\documentclass[oneside,reqno]{amsart}
\usepackage{palatino}
\usepackage[T1]{fontenc}
\usepackage[latin1]{inputenc}
\usepackage{amssymb}
\IfFileExists{url.sty}{\usepackage{url}}
                      {\newcommand{\url}{\texttt}}

\makeatletter

 \theoremstyle{plain}
 \theoremstyle{definition}
 \newtheorem*{defn*}{Definition}
 \theoremstyle{plain}
 \newtheorem{thm}{Theorem}
 \newtheorem*{claim}{Claim}
 \theoremstyle{plain}
 \newtheorem{lem}{Lemma}
 \theoremstyle{plain}
 \newtheorem*{conjecture*}{Conjecture}
 \theoremstyle{remark}
 
 \theoremstyle{definition}

\newtheorem{cor}{Corollary}
 
 \newcommand{\Z}{\ensuremath{\mathbb{Z}}}
 \newcommand{\N}{\ensuremath{\mathbb{N}}}
 
 \newcommand{\E}{\ensuremath{\mathbb{E}}}
 \newcommand{\PP}{\ensuremath{\mathbb{P}}}
\DeclareMathOperator{\Wait}{Wait}
\DeclareMathOperator{\Vis}{Vis}

\ifx\hyperlink\undefined
\newcommand{\hrlb}{ }
\newcommand{\ecp}[1]{ \url{http://www.math.washington.edu/~ejpecp/ECP/viewarticle.php?id=#1&layout=abstract}}
\newcommand{\ejp}[1]{ \url{http://www.math.washington.edu/~ejpecp/viewarticle.php?id=#1&layout=abstract}}
\else
\newcommand{\hrlb}{\\}
\newcommand{\ecp}[1]{\\ 
\href{http://www.math.washington.edu/~ejpecp/ECP/viewarticle.php?id=#1&layout=abstract}{\nolinkurl{http://www.math.washington.edu/~ejpecp/ECP/index.php}}}
\newcommand{\ejp}[1]{\\ 
\href{http://www.math.washington.edu/~ejpecp/viewarticle.php?id=#1&layout=abstract}{\nolinkurl{http://www.math.washington.edu/~ejpecp/index.php}}}
\fi
\newcommand{\mydot}{%
\begin{picture}(1,1)(0,0)
\put(3,0.5){\circle*{1}}
\end{picture}
}
\makeatother
\begin{document}
\title{Excited random walk against a wall}

\begin{abstract}
We analyze random walk in the upper half of a three
dimensional lattice which goes down whenever it
encounters a new vertex, a.k.a.~excited random walk. We show that it
is recurrent with an expected number of returns of $\sqrt{\log t}$.
\end{abstract}
\author{Gideon Amir}
\address{GA: Weizmann Institute, Rehovot, 76100, Israel}
\email{gideon.amir@weizmann.ac.il}
\author{Itai Benjamini}
\address{IB: Weizmann Institute, Rehovot, 76100, Israel}
\email{itai.benjamini@weizmann.ac.il}
\author{Gady Kozma}
\address{GK: Weizmann Institute, Rehovot, 76100, Israel}
\email{gady\makebox[1.2 ex]{\mydot}kozma@weizmann.ac.il}
\maketitle
\section{Introduction}
The model we will analyze in this paper (see section
\ref{sec:prelim} for a precise definition) is a variation on
excited random walk. Excited random walk is a walk on a $d$-dimensional
lattice ($d=1$ seems to be the richest case) which has 
a drift in some fixed direction whenever it encounters a new vertex. See
\cite{BW03, V03, K, Z05, Z} for recent results, \cite{PW97, D99} for a
Brownian motion analog, and \cite{AR} for some simulation results.
Excited random walk is proving to be far more tractable than other
self interacting processes such as the reinforced random walk or
the ``true'' self-avoiding walk.

In this paper we shall perform excited random walk on a half
space. Thus the walk's natural drift downward is counterbalanced
by the stiff floor. In a sense, the walk exhibits a self critical
behavior: if the walk ``tries to escape'' and visits a large
number of new vertices, it is pushed down to the floor and becomes
almost 2 dimensional and thus recurrent. If on the other hand the
walk returns too many times to the same vertices, it will upon
returning typically behave like simple random walk in 3 dimensions, which
is transient.


Thus, if a two dimensional random walk has approximately $\log t$ returns to
the origin until time $t$, and a three dimensional random walk has
approximately $1$ such return, we should expect excited random walk to
take some intermediate value. A somewhat less vague, but still
heuristic argument, says that the value should be $\sqrt{\log t}$: the
projection of the walk on the $(x,y)$ plane is a two dimensional
random walk so it returns to every column about $\log t$
times. If it reaches $x$ vertices in the column, it would accumulate a
downward drift of $x$. Assuming homogeneousness, it would visit the
floor about $\log t/x$ times and accumulate this amount of upward
drift. Since these should balance we get $x=\log t/x$ or
$x=\sqrt{\log t}$.

We shall prove that the $\sqrt{\log t}$ heuristic is in fact
accurate, and get in particular that the walk is recurrent, a fact
which is not at all clear a-priori. The actual proof only follows
the heuristic half way. The proof of the upper bound (see section
\ref{sec:upper}) will use different methods. The proof of the
lower bound will mimic the heuristic argument, but will use the
already established upper bound. An important tool in proving the
lower bound will be a coupling argument between two instances of
excited random walk (section
\ref{se:coupling}), which enables us to strengthen the upper bound
and replace the role of homogeneousness in the heuristic.
Unfortunately, we were forced to assume a deterministic downward
drift to make ends meet.

As a side remark, the square root heuristic also works for the
analogous model in two
dimensions, and one gets that the average number of visits of excited
random walk to the floor is of the order of $\sqrt[4]{t}$. The
two dimensional case is less interesting because recurrence
can be proved easily by coupling to simple random walk (so that the
simple random walk is always higher than the excited random walk) and
this argument does not
require deterministic drift or specific floor behavior. We will not present
any details of the two dimensional case.



\subsection{Open Problems}
As already remarked, we were not able to prove the case where the
walk, upon hitting a new vertex, goes down with some probability
$p<1$. The upper bound (theorem
\ref{thm:upper}) carries through unchanged, but the coupling argument
(lemma \ref{lem:cpl}) requires that the configuration be downward
closed which is not true for a probabilistic drift meaning that our
techniques only give a lower bound for the \emph{average} number of visits to
the floor, but one cannot deduce from that a lower bound for the
number of visits to a specific vertex, or even recurrence.

\subsection{Acknowledgments}
We would like to thank Alain-Sol Sznitman for a useful discussion
concerning the heuristic argument above. This research was carried
out while GK was staying at the Weizmann Institute of Science as a
Charles Clore postdoctoral fellow, and he would like to thank the
institute for its generous support. The writing process was partially
supported by the National Science Foundation under agreement
DMS-0111298. Any opinions, findings and conclusions or recommendations
expressed in this material are those of the authors and do not
necessarily reflect the views of the National Science Foundation.

\section{Preliminaries}\label{sec:prelim}

\begin{defn*}
In this paper, excited random walk (ERW) is a process of points $\{
R(t)\linebreak[1]=(x(t),y(t),z(t))\}_{t=1}^{\infty}$ in
$\mathbb{H}:=\{(x,y,z)\in\mathbb{Z}^{3}:z\geq0\}$
such that $R(0)=\underline{0}=(0,0,0)$ and $R(t+1)$ is created as
follows:
\begin{itemize}
\item Floor --- when the walk is currently on a floor vertex,
i.e.~$z(t)=0$ it moves with probability $\frac{1}{5}$ up and with
probability $\frac{1}{5}$ to either of the 4 sides.
\item Visited --- when the walk is on a non-floor vertex it visited
before, namely $R(t)=R(u)$ for some $u<t$, then it moves like a simple
random walk.
\item New --- when the walk is at a non-floor vertex it never visited,
then it moves downwards deterministically, namely
$R(t+1)=(x(t),y(t),z(t)-1)$.
\end{itemize}
We can also talk about an ``ERW starting from $v$'' for some
$v\in\mathbb{H}$ and in this case we take $R(0)=v$ instead.
\end{defn*}

A more-or-less equivalent process is the symmetric ERW, defined
on $\mathbb{Z}^{3}$ with the {}``excited steps'', i.e.~the steps
performed when reaching a new vertex $(x,y,z)$ go down if $z>0$
and up if $z<0$. Also, a vertex $(x,y,z)$ is considered visited if
$R(u)=(x,y,z)$ or if $R(u)=(x,y,-z)$ for some $u<t$. At the middle level
the walk has probability $\frac{1}{5}$ for the sides and $\frac{1}{10}$
for the up or down. Thus
if $R(t)$ is a symmetric ERW then
$(x(t),y(t),\left|z(t)\right|)$
is an ERW, and vice versa, an ERW can be symmetrized
by adding random coin flips that will decide, whenever the walk is
at $z=0$ whether to go up or down.

Given an ERW $R$ at some time $t$, we denote by
$\Vis_R(t)$ the set of visited (non-floor) vertices
i.e.~$\{R(u):u<t\}\setminus \{(x,y,0)\}$. When $t$
is clear from the context, we shall omit it, referring to
the set of visited vertices as $\Vis_R$.

It is important to notice that conditioning on the past (i.e.\ on $R[0,t]$ for
some $t$) is identical to conditioning on $\Vis(t)$ and $R(t)$. This can be
given a formal meaning using the notion of conditional probability (see e.g.\
\cite[4.1c]{D96}) but instead we shall use the following definition.
For any set of vertices $\mathcal{V}\subset\mathbb{H}$
we shall define an ERW ``starting from $(v,\mathcal{V})$'' by defining
$R(0)=v$ and $\Vis_R(0):=\mathcal{V}$ and continuing in the natural
way. Clearly, $R[t,\infty)$ is the same as an ERW starting from
$(R(t),\Vis_R(t))$. Some of our results (mainly theorem \ref{thm:upper}) hold
also for an ERW starting from $\mathcal{V}$ for a general $\mathcal{V}$. In
others (like theorem \ref{thm:lower}) this generalization requires
assuming that $\mathcal{V}$ is a legal configuration, i.e.\ a configuration
which can be $\Vis_R(t)$ with positive probability. It is easy to see that
this is equivalent to $\mathcal{V}$ being finite, connected and
downward-closed.

\medskip

The \emph{hitting time} of a random walk $R$ (simple or excited) of a subset
$A\subset \Z^3$ is the first time $t>0$ such that $R(t)\in A$. Notice
that due to the requirement $t>0$ the hitting time is non-trivial even
if $A$ contains the starting point $R(0)$.

For a subset $A\subset \Z^3$ we denote by $\partial A$ the internal
boundary of $A$ i.e.~all vertices in $A$ with a neighbor outside $A$.
$B(v,r)$ will denote a ball around $v$ with radius $r$.

When we write $\log n$ we always mean $\max\{1,\log n\}$ and $\log
0:=1$. We use $C$ and $c$ to denote various universal positive constants, 
which could take different values, even inside the same formula.
$C$ will be used for constants which are {}``large enough'' and $c$
for constants which are {}``small enough''.

\section{Upper Bound on the number of visits to a point}
\label{sec:upper}

In this section we shall prove the following:

\begin{thm}\label{thm:upper}
Let $R$ be an ERW starting from some point $w\in\mathbb{H}$ and some
configuration $\mathcal{V}$. Let $v\in\mathbb{H}$. Let $V(t)=V(t;v)$
be the number of times $R$ returns to $v$ until time $t$. Then\[
\mathbb{E}V(t)\leq C\sqrt{\log t}.\]

\end{thm}
We shall use the theorem for $\mathcal{V}$ empty or the visited set of
some past ERW. It should be noted, however, that the theorem holds
when starting from any configuration, even one which is impossible to
achieve using ERW, for example an isolated visited vertex.

The proof requires dividing into subshells and doing careful analysis of the
transitions of ERW from one shell to the next. This is somewhat technical,
so we shall first (section \ref{sub:simp})
sketch a ``single shell'' (well, a ball) argument which gives a weaker
result with some simplifying assumptions. We hope this makes the proof
clearer. Next (section \ref{sub:oned}) we shall give some simple lemmas that
allow to compare two one-dimensional processes. After that (section
\ref{sub:shells}) we shall analyze ERW in a single shell, and finally (section
\ref{sub:proof1}) we shall wrap the whole thing up.

We shall prove the claim for a symmetric ERW, and lose only
a factor of two in the constant $C$. Obviously, we may assume $t$
is large enough. It will be convenient to assume $V(t;w)$
also counts time $0$.

\subsection{Simplified sketch}
\label{sub:simp}

The argument we will now sketch only gives that the number of returns
to a vertex
$v$ is $\leq C\log^{5/6}t$. To see this, examine symmetric ERW and
take a ball around $v$
with radius $\log^{1/3}t$ and examine ``visits to the ball'' by which
we mean the time between one hitting of the ball and the next escaping
from a ball of double radius. Inside the ball there are only $\log
t$ vertices therefore, if we think about visits to the ball that
pass through more than $\log^{1/6}t$ new vertices as {}``bad'',
then there can be no more than $\log^{5/6}t$ bad visits. On
the other hand, a {}``good'' visit is limited by the amount it
differs from a simple random walk: even if we allow an opponent
with total view of past and future to distort a simple random walk
at less than $\log^{1/6}t$ times, she cannot force the walk to
pass through $v$ unless the original, simple random walk, passes
through a ball of radius $\log^{1/6}t$ around $v$, the probability
for which (at every visit to the outer ball) is $\leq
C\log^{-1/6}t$. Since, from two dimensional arguments, there are
only $\log t$ such visits, then the total number of good visits to
$v$ is also bounded by $C\log^{5/6}t$. Summing the good and the bad
proves the claim.

\subsection{One dimensional processes}
\label{sub:oned}

\begin{lem}
\label{lem:obvious} Let $R$ be a nearest-neighbor stochastic process
on $\{1,\dotsc,n\}$ with a uniform bound
\[
\PP\big(R(t)=R(t)+1\,\big|\,R[0,t]\big)\leq q_{R(t)}\quad 
\forall R(t)\in\{2,\dotsc,n-1\}
\]
for some numbers $\{q_2,\dotsc,q_{n-1}\}$.
Let $S$ be a nearest-neighbor Markov chain on $\{1,\dotsc,n\}$
starting from $R(0)$ with the transition probabilities (again from $i$
to $i+1$) equal to $q_i$. Let $s$ and $r$
be the probabilities that $S$ and $R$ respectively hit $n$ before
hitting $1$. Then
\[
r\leq s.
\]
\end{lem}

\begin{proof}
This follows directly from the fact that we can couple the two
processes in a way such that $S$ is always to the right of $R$ and
their difference is always even. See 
e.g.~\cite{L02} for some background on the coupling method.
\end{proof}

\begin{lem}
\label{lem:exp}Assume $q_{i}\leq q<\frac{1}{2}$, $i=2,\dotsc,n-1$ are
the transition probabilities of a nearest-neighbor Markov chain $R$ on
$\{1,\dotsc,n\}$ and let $r_{j}$ be the probability that $R$ hits
$n$ before hitting $1$ if $R$ starts at $j$. Then
\[
r_{j+1}\geq r_{j}(1+c)
\]
where the constant $c>0$ may depend on $q$ but not on $n$.
\end{lem}

\begin{proof}
Let $S$ be a nearest-neighbor Markov chain on $\{1,\dotsc,j+1\}$
starting from $j$ with transition probabilities $q$. Then the
probability that $S$ reaches $1$ before $j+1$ can be calculated
explicitly \cite[example 5.3.5]{D96} and is $>(1-2q)/(1-q)$ which we may denote by $c$. Hence
(using lemma \ref{lem:obvious}) we get that the probability $p$ of $R$
starting from $j$ to reach $1$ before $j+1$ satisfies $p>c$. However,
\[
r_{j}=(1-p)r_{j+1}\]
and we are done.
\end{proof}

\begin{lem}
\label{lem:perturb}
Assume $0<p\leq q_{i}\leq q<\frac{1}{2}$ are the
transition probabilities of a nearest-neighbor Markov chain $R$ on
$\{1,\dotsc,n\}$ starting from $j$ and let $r_{j}$ be the probability
that $R$ hits $n$ before hitting $1$. Then
\[
r_{j}\leq C\left(\frac{p}{1-p}\right)^{n-j}\prod_{i=1}^{n}(1+C(q_{i}-p))
\]
where $C$ may depend on $p$ and on $q$ but not on $n$ or on the $q_i$-s.
\end{lem}

\begin{proof}
In the case $q_{i}=p$ for all $i$, the numbers $r_{j}$ satisfy a
simple quadratic recursion, namely
$r_{j+2}-\frac{1}{p}r_{j+1}+\frac{1-p}{p}r_{j}=0$, which can be
solved explicitly to show that $r_{j}\leq C(p/(1-p))^{n-j}$. Hence
it is enough to measure the effect of a change in one $q_{i}$
namely, to show that if $q_{i}'=q_{i}$ except for one $I$, and
$q_{I}'>q_{I}$ then \begin{equation} r_{j}'\leq
r_{j}(1+C(q_{I}'-q_{I})).\label{eq:oneqIqIp}\end{equation} Let
$A$ and $B$ be some parameters. The numbers \[
s_{j}:=\begin{cases}
r_{j} & j<I\\
Ar_{j}+B & j\geq I\end{cases}\]
are {}``harmonic'' except possibly at $I-1$ and $I$ (meaning that
$s_{j}=q'_{j}s_{j+1}+(1-q'_{j})s_{j-1}$)
hence for the values of $A$ and $B$ satisfying that $s_{I}=r_{I}$
and $s_{I}=q_{I}'s_{I+1}+(1-q_{I}')s_{I-1}$ we would get that $s_{i}$
is harmonic and as a consequence, $r_{j}'\equiv s_{j}/s_{n}$. These conditions give that for the 
required $s_{I+1}$, 
\[|s_{I+1}-r_{I+1}|\leq r_I\left|\frac{1}{q'_I}-\frac{1}{q_I}\right|+r_{I-1}\left|\frac{1}{1-q'_I}-\frac{1}{1-q_I}\right|
\]
and by lemma \ref{lem:exp} this gives $|s_{I+1}-r_{I+1}|\leq
Cr_I(q_{I}'-q_{I})$.
Hence ($A$ and $B$ are linear in $s_{I}$ and $s_{I+1}$) we get that
\begin{align*}
A & =1+O\left(q_{I}'-q_{I}\right)\frac{r_{I}}{r_{I+1}-r_{I}} & B &
=O(q_{I}'-q_{I})\frac{r_{I}^{2}}{r_{I+1}-r_{I}}\end{align*}
(the constant implicit in the notation $O$ here may also depend on
$p$ and $q$). Another appeal to lemma \ref{lem:exp} allows to replace
the $r_{I}/(r_{I+1}-r_{I})$ factors with a constant and we get that
$s_{n}=1+O(q_{I}'-q_{I})$. Hence we get (\ref{eq:oneqIqIp}) and
the lemma.
\end{proof}
\subsection{Behavior in shells}
\label{sub:shells}

The proof of theorem \ref{thm:upper} in the next section uses lemma
\ref{lem:perturb} for one-dimensional processes 
created by examining the hitting times of shells of radius $4^n$.
A three dimensional Brownian motion (which is our model, in some vague sense)
starting from a point on a shell of radius 
$r$ has probability $\frac{1}{5}$ to reach $\frac{1}{4}r$ before
$4r$, independently of the starting point. For an ERW this probability
depends on the starting point and on $\Vis$, but
we will compare the probability that a \emph{good visit}, in the same sense as
in the proof sketch above, reaches
$\frac{1}{4}r$ before $4r$ to $\frac{1}{5}$. This will be done in 
lemma \ref{lem:r4r}. The other lemmas handle boundary cases: lemmas
\ref{lem:cyl} and \ref{lem:Q} are for the outermost shell and lemma
\ref{lem:r1} is for the innermost shell.

\begin{lem}
\label{lem:cyl}Let $r\geq 1$ and let
$K_{r}:=B((x,y),r)\times\mathbb{Z}$
namely an infinite vertical cylinder around $v=(x,y,z)$. Define stopping times
$t_{i}^{\mathrm{in}/\mathrm{out}}$
as follows: $t_{0}^{\mathrm{out}}:=1$, and for $i\geq1$,
\begin{align}
  t_{i}^{\mathrm{in}} & :=\min\{ t\geq t_{i-1}^{\mathrm{out}}:R(t)\in
  K_{r}\}\label{eq:defti}\\
  t_{i}^{\mathrm{out}} & :=\min\{ t>t_{i}^{\mathrm{in}}:R(t)\not\in
  K_{2r}\}.\nonumber 
\end{align}
Define $l:=\max\{ i:t_{i}^{\mathrm{in}}<t\}$, {}``the number of visits
to $K$''. Then
\[
\mathbb{P}(l>\lambda\log t)\leq Ce^{-c\lambda}\quad\forall\lambda>0.
\]
\end{lem}
\begin{proof}
For every $i>1$ we know that $t_{i-1}^\textrm{out}$ is in the exterior
boundary of $K_{2r}$ which we denote by $S$. Hence for $i>1$
\[
\PP(t_{i}^{\textrm{in}}>t_{i-1}^{\textrm{out}}+t)\geq
\min_{w\in S}\PP^w(\textrm{$R$ hits $\partial K_{t+2r}$ before $K_r$})
\]
which is a purely two-dimensional question. Denoting by $a$ the
two-dimensional discrete harmonic potential we can continue the inequality 
(see \cite{K87} for a nice
exposition of the connection between harmonic functions and the
harmonic potential in particular and hitting probabilities)
\begin{align*}
&\geq \frac{\min\{a(w):w\in S\}-\max\{a(w):w\in \partial K_r\}}
{\max\{a(w):w\in \partial K_{t+2r}\}-\min\{a(w):w\in \partial K_r\}}\\
\intertext{and since $a(w)=A\log |w|+B+o(1)$ \cite[P12.3]{Sp}}
&=
\frac{\log 2+o(1)}{\log((t+2r)/r)+o(1)}\geq\frac{c}{\log
t}
\end{align*}
which gives the lemma immediately.
\end{proof}

\begin{lem}\label{lem:r4r}
Let $r>4$, let $v\in \mathbb{Z}^3$ and let $w\in \partial B(v,r)$.
Let $R$ be a (symmetric) ERW starting from $w$ and some configuration
$\mathcal{V}\subset\mathbb{Z}^3$, and let $T$ be its hitting time of
$\partial B(v,4r)\cup \partial B(v,\frac{1}{4}r)$. Let
$\epsilon>r^{-1/2}$ be some parameter. Let $G$ be the event that $R$
encounters less than $\epsilon r$ new sites until $T$. Then
\begin{equation}
\mathbb{P}(\{R(T)\in\partial B(v,\textstyle{\frac{1}{4}}r)\}\cap G)\leq
\textstyle{\frac{1}{5}}+C\epsilon^{1/4}.
\end{equation}
\end{lem}

\begin{proof}
Couple $R$ to a simple random walk $W$ in the following manner:
if $R$ is in a visited (non-floor) vertex let $R$ and $W$ perform the
same step. Otherwise they walk independently according to their
respective rules. The lemma will be mostly proved once we estimate
$|R(t)-W(t)|$. There are two sources for the discrepancy: new vertices
and floor vertices. Therefore let us write $|R(t)-W(t)|=|N(t)+F(t)|$ where
\[
N(t)=\sum_{\substack{u<t\\ R(u)\textrm{ is a new vertex}}} 
  R(u+1)-R(u)-W(u+1)+W(u)
\]
and $F$ is the same for floor vertices. Now, $G$ obviously implies
$|N(t)|\leq 2\epsilon r$ for all $t\leq T$ so we need only estimate
$F$. Now, for every time $t$ when $R$ is in a floor vertex, the
expected motion of $R$ is zero (remember that we are talking about
the symmetric ERW). In other words, if we denote by $t_i$ the $i$'th
hitting of the floor then $F(t_i)$ is a symmetric random walk on
$\mathbb{Z}^3$ with bounded steps. By the reflection principle
(see e.g.~\cite[chapter 2, lemma 1]{K85}) $\max_{i\leq n}|F(t_i)|$ has the
same tail 
behavior as $|F(t_n)|$ i.e.\ a square-exponential one. Denoting by $f$
the number of times $R$ hits the floor by time $T$ we get
\begin{equation}
\mathbb{P}(
  \{\max_{t\leq T} |F(t)|>\lambda \sqrt{N}\} \cap
  \{f\leq N\})
\leq
Ce^{-c\lambda^2} \quad
\forall\lambda>0,\,\forall N\in\mathbb{N}.\label{eq:Ff}
\end{equation}
Examine one time $t$ when $R(t)$ is at the floor. It is easy to see
that a simple random walk starting from a floor point has probability
$\geq c/r$ to exit the ball $B(v,4r)$ before returning to the floor.
Therefore we have
\[
\mathbb{P}(\textrm{$R$ hits $\partial B(v,4r)$ or some new vertex
  before returning to the floor})
\geq \frac{c_1}{r}.
\]
This implies that if $f>(2/c_1)\epsilon r^2$ 
then with probability $>1 - Ce^{-c\epsilon r}$ there are at least
$\epsilon r+1$ times $t_i<T$ when the event above happened. This,
however, contradicts the event $G$ so we get
\begin{equation}
\mathbb{P}(\{f>(2/c_1)\epsilon r^2\}\cap G)\leq Ce^{-c\epsilon r}.\label{eq:fr}
\end{equation}
Combining (\ref{eq:Ff}) and (\ref{eq:fr}) we get
\[
\mathbb{P}(\{
  \max_{t\leq T}|F(t)|>C\lambda r\sqrt{\epsilon}\}
  \cap G)\leq 
Ce^{-c\lambda^2}+Ce^{-c\epsilon r}
\quad\forall\lambda>0.
\]
We pick $\lambda:=\frac{1}{C}\epsilon^{-1/4}$ and get (using also the requirement $\epsilon > r^{-1/2}$)
\[
\mathbb{P}(\{\max_{t\leq T}|F(t)|>r\epsilon^{1/4}\}\cap G)\leq C\epsilon.
\]
This is the estimate of $F$ that we need.

Now, in general we have for any $s\geq 1$ that $R(T)\in \partial
B(v,\frac{1}{4}r)$ implies that either for some $t\leq T$ we have
$|N(t)+F(t)|>s$ or that $W$ hits $\partial B(v,\frac{1}{4}r+s)$ before
hitting $\partial B(v,4r+s)$. The probability for that to happen
(denote it by $q$) we calculate using the discrete Green function of $\Z^3$
(denote it by $a$) the same way we used the harmonic potential in the
previous lemma. Since $a(z)=\frac{c}{|z|}+O\left(\frac{1}{|z|^{3}}\right)$
\cite[Theorem 4.3.1]{L} we get $q=\frac{1}{5}+O(s/r)$.
Applying this with
$s=r(2\epsilon + \epsilon^{1/4})$ we get
\begin{align*}
\mathbb{P}(\{R(T)&\in \partial B(v,{\textstyle\frac{1}{4}}r)\}\cap G) \leq\\
&\leq \mathbb{P}(W\textrm{ hits }\partial B(v,{\textstyle\frac{1}{4}}r+s)) +
  \mathbb{P}(\{\max_{t\leq T}|N(t)+F(t)|>s\}\cap G) \leq\\
&\leq \frac{1}{5}+O(s/r)+C\epsilon \leq \frac{1}{5}+C\epsilon^{1/4}.
\qedhere
\end{align*}
\end{proof}

\begin{lem}
\label{lem:Q}
Let $r>4$, let $v=(x,y,z)$ and let $w\in \partial B(v,r)$.
Let $R$ be a (symmetric) ERW starting from $w$ and some configuration
$\mathcal{V}\subset\mathbb{Z}^3$, and let $T$ be its hitting time on
$\partial B(v,4r)\cup \partial B(v,\frac{1}{4}r)$. Denote 
\[
Q:=(\partial B(v,4r) \cap B((x,y),2r)\times \Z) 
  \cup \partial B(v,\tfrac{1}{4}r)
\]
i.e.~$Q$ is the union of a) the outer sphere intersected
with a concentric vertical cylinder of half its radius and b) the inner
sphere. Let $G$ 
be the event that $R$ encounters less than $\epsilon r$ new sites until $T$
for some $\epsilon$ sufficiently small. Then
\begin{equation*}
\mathbb{P}(\{R(T)\in Q\}\cap G)\leq 1-c
\end{equation*}
\end{lem}

The proof is very similar to the proof of the previous lemma
--- in fact, simpler --- the only additional fact needed is that a
\emph{simple} random walk has probability $<1 - c$ to hit $Q$. This is
quite easy to see and we omit any further details about the proof of
lemma \ref{lem:Q}.

\begin{lem}\label{lem:r1}
Let $r>1$, let $v\in \mathbb{Z}^3$ and let $w\in \partial B(v,r)$.
 Let $R$ be a (symmetric) ERW starting from $w$ and some configuration
 $\mathcal{V}\subset\mathbb{Z}^3$, and let $T$ be its hitting time of
 $\partial B(v,4r)\cup \{v\}$. Let $G$ be the event that $R$
 encounters no new sites until $T$. Then
\begin{equation}
\mathbb{P}(\{R(T)=v\}\cap G)\leq C/r.
\end{equation}
Similarly if $R$ starts from $v$ then this probability is $\leq 1-c$.
\end{lem}
\begin{proof}
This time we couple $R$ to a random walk $W$ which has the same
behavior as $R$ 
at the floor, i.e.~when $W$ hits the floor it has probability $\frac{1}{5}$ to
go to each of its floor neighbors, and probability $\frac{1}{10}$ for each of
its vertical neighbors, but other than that is simple. Clearly if $G$ happened
then $R(t)=W(t)$ for all $t$ so it is enough to estimate the corresponding
probabilities for $W$.

However $W$ is a reversible random walk (meaning that it can be
realized as a walk on a weighted graph) so Varopoulos \cite{V85} and
Hebisch and Saloff-Coste \cite[theorem 2.1]{HSC93} apply. Together
they give that the probability $p_t(x,y)$ that $W$ starting from $x$
and going $t$ steps will be at $y$ satisfies
\[
p_t(x,y)\leq \frac{C}{t^{3/2}}e^{-c|x-y|^2/t}\quad\forall x,y,t.
\]
Summing over $t$ gives that the discrete Green function
$a(x,y)=\sum_t p_t(x,y)$ satisfies $a(x,y)\leq C/|x-y|$ and
$a(x,x)\leq C$. $a$ is harmonic so the same calculations as in 
lemma \ref{lem:r4r} give the estimates for the probabilities.
\end{proof}

\subsection{Proof of theorem \ref{thm:upper}}
\label{sub:proof1}

Let $\beta_{j}=4^{-j}\sqrt{\log t}$ defined for $j=1,\dotsc,J$
until the first $J$ such that $\beta_{J}<\sqrt[6]{\log t}$, and
let $\beta_{J+1}=\beta_{J+2}=1$. Let $S_{i}:=\partial B(v,\beta_{i})$, and in
particular $S_{J+1}=S_{J+2}=\{ v\}$. The spheres $S_{i}$ are the analogue
of the sphere at $\log^{1/3}t$ discussed in the ``simplified sketch''
section. Let $t_{i}$ denote stopping times at these
spheres defined, somewhat similarly to (\ref{eq:defti}), by
\begin{align*}
t_{1} & =\min u:R(u)\in\bigcup_{j=1}^{J+1}S_{j}\\
t_{i+1} & =\min u>t_{i}:\begin{cases}
u\in S_{j-1}\cup S_{j+1} & \textrm{when}\quad R(t_{i})\in S_{j},2\leq j\leq J+1\\
u\in S_{2} & \textrm{when}\quad R(t_{i})\in S_{1}
\end{cases}\end{align*} (notice
the asymmetry at $v$ --- the only case where $R(t_{i})$ and
$R(t_{i+1})$ may belong to the same $S_{j}$). Let $\delta\in (0,1)$ be
some parameter sufficiently small to be fixed
later. In fact, it is enough to take $\delta:=\min\{(2/15C_{\text{lemma
      \ref{lem:r4r}}})^4,\epsilon_{\text{lemma \ref{lem:Q}}},\frac{1}{2}\}$,
  but the only meaning of this expression is in the various conditions that
  will appear below.
Let $G_{i}$ ($G$ standing for {}``good'') be the event that $R$ hits
less than $\delta j^{-8}\beta_{j+1}$ new sites between time $t_{i}$
and time $t_{i+1}$ where $j\leq J+1$ is given by $R(t_{i})\in
S_{j}$. Obviously, there is nothing stopping $\delta
j^{-8}\beta_{j+1}$ to be smaller than $1$ (indeed it must be if $j\geq
J$). Remembering the {}``simplified sketch'' section, the event
$G_{i}$ is the analogue of the event {}``good visit to $v$'' with
respect to the relevant sphere $S_{i}$.

To estimate $V(t)$ we examine the walk performed before the time
when $v$ was hit, and ask: when has $G_{i}^{c}$ (the complement
of $G_{i}$) occurred last? More precisely, define the event $H_{i}^{0}$
to be $G_{i}\cap G_{i+1}\cap\dotsc\cap G_{k-1}$ where $k\geq i$ is the
first such that $R(t_{k})\in S_{1}\cup\{ v\}$ and let $H_{i}$ be
the event that $H_{i}^{0}$ happened and $R(t_{k})=v$. In particular,
if $R(t_i)=v$ then $H_i$ happens while if $R(t_i)\in S_1$ it does
not. Define
\begin{enumerate}
\item $E_{j}$ ($2\leq j\leq J+1$) to be the number of $i$-s
  satisfying that $R(t_{i})\in S_{j}$, that $G_{i}$ did \emph{not}
  happen, and that $H_{i+1}$ did happen. 
\item $E_{J+2}$ to be the number of $i$-s such that $R(t_{i})=v$ and
  $G_i\cap H_{i+1}$ happened.
\item $E_{1}$ to be the number of $i$-s such that $R(t_{i})\in S_{1}$
  and $H_{i+1}$ happened.
\end{enumerate}
For all these we count only $i$-s that the relevant $t_{k}$, i.e.~the time
where $R(t_{k})=v$, happened before time $t$. Clearly,
$V(t)=\sum_{j=1}^{J+2}E_{j}$ therefore it is enough to estimate these
$E_{j}$-s. Now if $t$ is sufficiently large then we can apply lemma
\ref{lem:r4r} for all $j<J$ (the problem is only in the condition
``$\epsilon>r^{-1/2}$'' of lemma \ref{lem:r4r} where here
$\epsilon=\frac{1}{4}\delta j^{-8}$ and $r=\beta_j$). For $j=J,J+1$
we use lemma \ref{lem:r1} and in total we get
\begin{equation}
\mathbb{P}(G_{i}\cap\{R(t_{i+1})\in S_{j+1}\}\,|\,
R[0,t_{i}],R(t_{i})\in S_{j})\leq
\begin{cases}
  \frac{1}{5}+C\delta^{1/4}j^{-2} & 1\leq j< J\\
  C\log^{-1/6}t & j=J\\
  1-c & j=J+1.
\end{cases}
\end{equation}
Here and below we use the notation $\PP(X|Y,E)$ for a variable $Y$ and an event
$E$ to mean the function $\PP(X|Y)$ restricted to $E$ --- here everything is
discrete so this simply means that the inequality holds for any value of
$R[0,t_i]$ for which $R(t_i)\in S_j$.
Denote the values on the right hand side by $q_j$. This allows us to
estimate $\mathbb{P}(H_{i}\,|R[0,t_{i}])$ 
by comparing the process $R(t_{i})$ to a Markov chain $Q_{j}$ on
$\{1,\dotsc,J+1\}$ starting from $j$ with the transition probabilities
$q_{j}$ (we use here lemma \ref{lem:obvious}).
If $\delta$ is sufficiently small (explicitly if 
$\delta\leq(2/15C_{\text{lemma \ref{lem:r4r}}})^4$) we would have 
$q_{j}\leq\frac{1}{3}$ for all $j<J$. Hence we can use lemma \ref{lem:perturb}
on the interval $[1,J]$ and we get

\[
\mathbb{P}(Q_{j}\textrm{ hits }J\textrm{ before }1)\leq
C4^{j-J}\prod_{j=1}^{J}(1+Cj^{-2})\leq
C4^{j-J}=C\frac{\beta_{J}}{\beta_{j}}.
\]
The step from $J$ to $J+1$ contributes another $C\beta_{J}^{-1}$
factor so we end up with
\begin{alignat}{2}
\mathbb{P}(H_{i}\,|\, R[0,t_{i}]) &
\leq\frac{C}{\beta_{j}} &\quad&
\textrm{when }R(t_{i})\in S_{j},j\leq J\label{eq:Hismall}\\
\mathbb{P}(G_i\cap H_{i+1}\,|\, R[0,t_{i}]) &
\leq 1-c &&
\textrm{when }R(t_i)=v.\label{eq:Hivsmall}
\end{alignat}

We note that $B(v,\beta_{j-1})\setminus B(v,\beta_{j+1})$ contains
$<C\beta_{j-1}^{3}$ points. Therefore the number of $i$-s such that
$G_{i}^{c}$ can occur together with $R(t_i)\in S_j$ is no more
than
\begin{alignat}{2}
\frac{C\beta_{j-1}^{3}}{\delta j^{-8}\beta_{j+1}} & \leq 
  C j^8 4^{-2j}\log t\quad & \textrm{when } j & <J\label{eq:volGc}\\
\beta_{j-1}^{3} & \leq C\sqrt{\log t} & \textrm{when } j & =J,J+1.\nonumber
\end{alignat}
(here and below we will be ``folding'' the $\delta$ into the
constants $c$, $C$). Using (\ref{eq:Hismall}) for $i+1$ (and $j+1$, which
would also estimate the case that $R(t_{i+1})\in S_{j-1}$) shows that
$E_j$ is dominated by a sum of independent Bernoulli trials with
probability $C\beta_{j+1}^{-1}$, so
\begin{alignat}{2}
\mathbb{P}(E_{j}>Cj^8 4^{-j}\sqrt{\log t}+
    \lambda j^4 2^{-j}\log^{1/4}t) &
  \leq Ce^{-c\lambda} \qquad &
  2 & \leq j<J\label{eq:Ej1Jm1}
\end{alignat}
while for $j=J,J+1$ we have deterministically
\begin{align}
E_{j}\leq C\sqrt{\log t}.
\label{eq:EjJJp1}
\end{align}
In
particular, $\mathbb{E}\sum_{j=2}^{J+1}E_{j}\leq C\sqrt{\log t}$.

Next we estimate $E_{1}$. We start with an estimate of the number
of $i$-s such that $R(t_{i})\in S_{2}$. Denote it by $F$. Define
$Q=(S_{1}\cap K_{2\beta_{2}})\cup S_{3}$ where $K$ is an infinite
cylinder as in lemma \ref{lem:cyl}. By lemma \ref{lem:Q}, if only
$\delta$ is sufficiently small 
($\delta\leq\epsilon_{\text{lemma \ref{lem:Q}}}$), 
\[
\mathbb{P}(G_{i}\cap\{R(t_{i+1})\in
Q\}\,|\, R[0,t_{i}],\,R(t_i)\in S_2)\leq 1-c
\]
Hence the
number $X_{1}$ of times this event happened satisfies
\begin{equation}
\mathbb{P}(X_{1}>(1-c)F+\lambda\sqrt{F\log F})\leq Ce^{-c\lambda^{2}}.
\label{eq:X1F}
\end{equation}
To prove (\ref{eq:X1F}) compare to an infinite sequence of 
Bernoulli trials $\epsilon_i$ with probability $1-c$ for
which a rough estimate (by summing over $s$) shows that
$\mathbb{P}(\exists s\leq t:\sum_{i=1}^s\epsilon_i > (1-c)s 
+\lambda \sqrt{s\log s})\leq Ce^{-c\lambda^2}$.

Next, the number of times $G_{i}^{c}$ happened is bounded using
(\ref{eq:volGc}) by $C\log t$ so we get that the number $X_{2}$ of
$i$-s for which $R(t_{i})\in S_{2}$ and $R(t_{i+1})\not\in Q$
satisfies 
\begin{equation} 
\mathbb{P}(X_{2}\leq cF-C\log t-\lambda\sqrt{F\log F})\leq
Ce^{-c\lambda^{2}}.\label{eq:X2small}
\end{equation} 
However, every
such event is an {}``entry into $K_{\beta_{2}}$'' in the sense of
lemma \ref{lem:cyl} so
\begin{equation}
\mathbb{P}(X_{2}>\lambda\log t)\leq Ce^{-c\lambda}\label{eq:X2big}.
\end{equation}
We get
\begin{align}
\mathbb{P}(F>\lambda\log t)&\leq
  \mathbb{P}(X_2>\sqrt{\lambda}\log t) + 
    \mathbb{P}(\{F>\lambda\log t\}\cap\{F>\sqrt{\lambda}X_2\}) \leq\nonumber\\
&\!\stackrel{(\ref{eq:X2big})}{\leq} Ce^{-c\sqrt{\lambda}} + 
    \mathbb{P}(\{F>\lambda\log t\}\cap\{F>\sqrt{\lambda}X_2\}).
\label{eq:Fsmall1}
\end{align}
For $\lambda$ sufficiently large we have that $F>\lambda\log t$ and
$F>\sqrt{\lambda}X_2$ imply that in fact $F>CX_2+C\log
t+C\sqrt{\lambda F \log F}$ and so by (\ref{eq:X2small})
\begin{equation}
  \mathbb{P}(\{F>\lambda\log t\}\cap\{F>\sqrt{\lambda}X_2\})\leq
  Ce^{-c\sqrt{\lambda}}.
\label{eq:Fsmall2}
\end{equation}
But (\ref{eq:Fsmall2}) can be made to hold not just for $\lambda$
sufficiently large by increasing the $C$ on the right hand side and
with (\ref{eq:Fsmall1}) we get
\begin{equation}
\mathbb{P}(F>\lambda\log t)\leq   Ce^{-c\sqrt{\lambda}}\quad\forall\lambda>0.
\label{eq:Fsmall}
\end{equation}
This is the estimate of $F$ that we need.

On the other hand, let $i$ satisfy the requirements for $E_{1}$,
namely $R(t_{i})\in S_{1}$ and $H_{i+1}$ has occurred. 
Using
(\ref{eq:Hismall}) we get\begin{align*}
\mathbb{P}(H_{i+1}\,|\, R[0,t_{i}],\, R(t_{i})\in S_{1})
 & =\mathbb{EP}(H_{i+1}\,|\, R[0,t_{i+1}],\, R(t_{i})\in S_{1})\leq\\
 &
\stackrel{(\ref{eq:Hismall})}{\leq}\mathbb{E}(C/\beta_{2})=C/\beta_{2}\end{align*}
where $\mathbb{E}$ here is the conditional expectation over the variable
$R[0,t_{i}]$. As in (\ref{eq:X1F}) above, we get
\begin{equation}
\mathbb{P}(E_{1}>CF/\beta_{2}+\lambda\sqrt{(F/\beta_{2})\log F})\leq
Ce^{-c\lambda}.\label{eq:E0F}\end{equation}
Combining (\ref{eq:Fsmall}) and (\ref{eq:E0F}) gives
\begin{multline}
\mathbb{P}(E_{1}>\lambda\sqrt{\log t}) \leq
  \mathbb{P}(\{E_1>\sqrt{\lambda}(F/\beta_{2})\}\cap
             \{E_1>\lambda\sqrt{\log t}\})\; +\\
+\; \mathbb{P}(F/\beta_{2}>\sqrt{\lambda\log t})\leq
Ce^{-c\lambda^{1/4}}.
\label{eq:E0}
\end{multline}
Hence $\mathbb{E}E_{1}\leq C\sqrt{\log t}$ and this part is estimated
as well.

The estimate of $E_{J+2}$ comes from (\ref{eq:Hivsmall}): again by
comparing to a sum of independent Bernoulli trials we get 
\[
\mathbb{P}(E_{J+2}\geq(1-c)V(t)+\lambda\sqrt{V(t)\log V(t)})
\leq Ce^{-c\lambda^{2}}
\]
or, equivalently,
\[
\mathbb{P}(E_{J+2}\geq C(V(t)-E_{J+2})+\lambda\sqrt{V(t)\log V(t)})\leq
Ce^{-c\lambda^{2}}.
\]
Adding (\ref{eq:Ej1Jm1}), (\ref{eq:EjJJp1}) and (\ref{eq:E0}) gives
\[
\mathbb{P}(V(t)-E_{J+2}\geq\lambda\sqrt{\log t})=
\mathbb{P}\Big(\sum_{j=1}^{J+1}E_j\geq\lambda\sqrt{\log t}\Big)
\leq Ce^{-c\lambda^{1/4}}\]
so
\begin{multline}
  \mathbb{P}(E_{J+2}\geq\lambda\sqrt{\log t}) \leq
    \mathbb{P}(\{E_{J+2}\geq\sqrt{\lambda}(V(t)-E_{J+2})\}\cap
    \{E_{J+2}\geq\lambda\sqrt{\log t}\})\;+\\
  +\;\mathbb{P}(V(t)-E_{J+2}\geq\sqrt{\lambda\log t})\leq
  Ce^{-c\lambda^{1/8}}\label{eq:EJp2}
\end{multline}
which shows that $\mathbb{E}E_{J+2}\leq C\sqrt{\log t}$ and since
this is the last term in $V(t)$, the theorem is proved.\qed

\begin{cor}\label{cor:cond}
For every $v\in\mathbb{H}$ we have
\begin{equation}\label{eq:Vncond}
\mathbb{E}(V(t;v)\,|\,V(t;v)\neq 0)\leq C\sqrt{\log t}
\end{equation}
where $V(t;v)$ is the number of visits to $v$ after $t$ steps.
\end{cor}
\begin{proof}
This is because conditioning on $V(t;v)\neq 0$ is identical to
considering an unconditioned ERW starting from $(v,\mathcal{V})$
where $\mathcal{V}:=\Vis(T)$, $T$ being the hitting time of $v$ and
walking for a distance of $t-T$. Applying theorem \ref{thm:upper}
shows that $\mathbb{E}V(t-T;v) \leq C\sqrt{\log t}$ for any $\mathcal{V}$ and
integrating over $T$ and $\mathcal{V}$ shows (\ref{eq:Vncond}).
\end{proof}
\begin{cor}[exponential decay of $V$]\label{exp_decay}
There exist constants $c,C$ s.t. for any point $v$ and any
$\lambda > 0$, $\PP(V(t;v)
> \lambda \sqrt{\log t}) < Ce^{-c\lambda}$.
\end{cor}
\begin{proof}
Using theorem \ref{thm:upper} and Markov's inequality we get some constant $K$
such that for every configuration $\mathcal{V}$ one has that an ERW starting
from $(v,\mathcal{V})$  has probability $<\frac{1}{2}$ to visit $v$ more than
$K\sqrt{\log t}$ visits in the next $t$ steps. Define $L:=\lfloor K\sqrt{\log
  t}\rfloor+1$ and let $\tau_k$ be the $kL$'th return to $0$ (here
$\lfloor\cdot\rfloor$ stands for the integer value). As in the previous corollary, the ERW after $\tau_k$ is the same as an ERW starting from $(v,\Vis(\tau_k))$ so we get
\[
\mathbb{P}(\tau_{k+1}>\tau_k+t\,|\,R[0,\tau_k])>\textstyle{\frac{1}{2}}.
\]
Hence we get that $\mathbb{P}(\tau_k<t)\leq 2^{-k}$, which was to be proved.
\end{proof}

\subsection{Postfix remarks}

The values chosen for the $\beta_{j}$ are in some sense
{}``non-opti\-mal''. A more natural choice would be
$\beta_{j}=e^{-2^{j}}\sqrt{\log t}$, i.e.~a doubly exponential
decreasing sequence. For example, if one decides to use only a
finite number of $\beta$-s (finite in the sense that the length
$J$ is independent of $t$) and looks for the optimal $\beta$-s,
the optimality requirement gives a set of equations which, when
solved, give a doubly exponential decreasing sequence with
$\beta_{1}=\log^{1/2-\epsilon}t$ and $\beta_{J}=\log^{1/6}t$.
Actually, the fact that we stopped our series $\beta_{j}$ when
reaching $\sqrt[6]{\log t}$ is an atavism from this optimization.
Either choice for the $\beta_{j}$ would give the same
conclusion in the theorem.

Lemmas \ref{lem:r4r}, \ref{lem:Q} and \ref{lem:r1} could have been simplified
significantly if the behavior of the ERW at the floor would have been
$\frac{1}{6}$ for its floor neighbors and $\frac{1}{3}$ for its upper
neighbor. Unfortunately, the coupling argument used in the next
section requires the probability of the upper neighbor to be 
$\leq\frac{1}{5}$.

Since corollary \ref{exp_decay} gives a very simple argument for the exponential decay of $V$, one might wonder why did we bother with all the intermediate estimates of the form
$\mathbb{P}(\textrm{something}>\lambda\sqrt{\log t})\leq
C\exp\left(-c\lambda^{\textrm{some fraction}}\right)$,
namely (\ref{eq:Ej1Jm1}),
(\ref{eq:EjJJp1}) or (\ref{eq:E0})? They seem to be necessary for
the calculation of $E_{J+2}$, (\ref{eq:EJp2}). We would like to
see a proof that can estimate $\mathbb{E}E_{J+2}$ using only
$\sum_{j=1}^{J+1}\mathbb{E}E_{j}$, but we were not able to
overcome some dependency issues.

\begin{conjecture*}
The correct tail decay is square-exponential, namely
$\mathbb{P}(V(t)>\lambda\sqrt{\log t})\leq Ce^{-c\lambda^{2}}$.
\end{conjecture*}
One possible interpretation of the word {}``correct'' above is:
for every $\lambda$ and every $t>t_0(\lambda)$,
$\mathbb{P}(V(t)>\lambda\sqrt{\log t})>ce^{-C\lambda^{2}}$.

\section{The coupling argument } \label{se:coupling}

As we will show below, when the starting configurations are downward closed,
it is possible to couple two instances of ERW such that one is always above
the other. Here it is more convenient to think about them as walks in a half
space rather than as the symmetrized version we used in the previous chapter,
so from now on we will use the half space version of ERW. The following lemma
uses this argument to show a certain monotonicity in the hitting
probabilities. It will be crucial towards the end.

\begin{lem}\label{lem:cpl}
Let $R,S$ be two ERWs, starting from a
$w\in \mathbb{H}$, and from visited configurations satisfying
$\Vis_R(0)\subset\Vis_S(0)$ which are both downward-closed. Let $v\in
\mathbb{H}$ be a floor vertex and let $V_R(t)$ and $V_S(t)$ be the
number of visits of $R$ and respectively $S$ to $v$ in the first $t$
steps. Then for any $t\in \N$ and  $k\in \N$ we have
$\PP(V_R(t) \geq k) \geq \PP(V_S(t) \geq k)$, and in
particular $\E(V_{R}(t)) \geq \E(V_{S}(t))$.
\end{lem}
\begin{proof}
We define a coupling between $R$ and $S$ so that for any instance
of the coupling the number of times $R$ hits $v$
before time $t$ is greater or equal to the number of times $S$
hits $v$ before time $t$. The coupling requires a time change so, if
we denote by $\tau$ the number of coupling steps, we need two time
change functions $t_R(\tau)$ and $t_S(\tau)$ to get back the time for
each process. For brevity, we will replace $R(t_R(\tau))$ with just
$R(\tau)$ or just with $R$ (ditto for $S$).

To define the coupling recall the three types of vertices an ERW can
be at: floor, visited and new. We define the coupling according to the types
of the vertices both walks are at, generally trying to make them walk
``together'':
\begin{itemize}
\item If both $R$ and $S$ are at the same type of vertex --- they move
together (i.e.~make the same step).
\item If one of them is at a new vertex, and the other is not --- the
one at the new vertex makes a move downwards, while the other one
waits.
\item If one of them is at a visited vertex, and the other at a floor
vertex, we let the first one move. If the move it made was downwards
--- the second walk waits. Otherwise, the second walk moves in the
same direction.
\end{itemize}

We denote by $\Wait_R(\tau)$ $(\Wait_S(\tau))$ the number of times
$R$ ($S$) waited until time $\tau$ of the coupling. Thus
$t_R(\tau)=\tau-\Wait_R(\tau)$ that is until
step $\tau$ of the coupling the walk $R$ makes $\tau-\Wait_R(\tau)$
real steps.

As above, when looking at a specific step $\tau$
of the coupled walk, we omit the step index from the various
values thus writing $x_R$ and $\Wait_R$ instead of $x_R(\tau)$ and
$\Wait_R(\tau)$.

The lemma will now follow from the following claim:

\begin{claim}\label{clm:cpl}
At each step of the coupling we have:
\begin{enumerate}
    \item  $x_R = x_S$ and $y_R = y_S$.
    \item  $z_R - \Wait_R = z_S - \Wait_S$.
    \item  $z_R \leq z_S$.
    \item  If $S$ is not at a new vertex, then $\Vis_S$ is downward
           closed and $\Vis_R\subset \Vis_S$.
\end{enumerate}
\end{claim}
\begin{proof}
First notice that both walks make the same moves on the
$(x,y)$-plane, regardless of the vertex type they are at, so at
each step $x_R = x_S$ and $y_R = y_S$, giving (1).

(2) follows from the fact that at each step of the coupling the
walks either move together, or one of them waits while the other
moves down.

To prove items (3) and (4) we use induction on the step of the
coupling. Assume that the claim holds up some step, and look at
the next step of the coupling.

We will first prove item (3) continues to hold. If both walks make
the same move, (3) continues to hold. Otherwise, we are in one of
the following two situations:
\begin{enumerate}
\item One of the walks is on a floor point, and the other one is
above it.
\item One of the walks is at a new vertex and the other is at a
visited one.
\end{enumerate}
In the first case, either the walks move together, or the walker at
the floor waits, while the second one goes down a step. But since
this means before the step $R$ was strictly below $S$, we get that
after the downward move still $z_R \leq z_S$. In the second case,
if $R$ is the one at the new vertex, or if $z_S > z_R$ then after
the next step still $z_R \leq z_S$, so the only case we must worry
about is that both walks are currently at the same vertex $q$, and
$q$ is a new vertex for $S$, while $R$ has already visited it. To
rule out this case, look at the first time $R$ visited $q$. At
that time, by the induction hypothesis, (items (1) and (3)), $S$
was directly above $R$, and by the coupling rules, $R$ would drop
at least one step, and  $S$ would drop until it reached a vertex
it has visited before (or the floor), prior to $R$ making any
sideways or upward move (and therefore prior to $R$ returning to $q$
and thus strictly before our current time). Thus when $S$ reaches
a non-new vertex, by the induction hypothesis (item (4)), $q\in
\Vis_R\subset \Vis_S$ contrary to our assumption. Thus
this last case is dismissed and we have proven (3).

To see (4), roll back to the last induction step $\sigma$ when $S$
is not in a new vertex. By the induction hypothesis,
$\Vis_S(\sigma)$ is downward closed and $\Vis_R(\sigma) \subset
\Vis_S(\sigma)$. Seeing that $\Vis_S$ remains downward closed is
obvious. To see $\Vis_R \subset \Vis_S$ divide into cases
according to which one moves. The only interesting case is when
they make a simultaneous move in the $x,y$ plane. However, if
$R(\sigma + 1)\not \in \Vis_S(\sigma)$ then so is $S(\sigma + 1)$
(because $\Vis$ is downwards closed) and then $S$ must drop until
closing $\Vis(S)$ before $\tau$, so $\Vis_S$ contains the entire
column $[0,z_S(\sigma + 1)]$ above $x_S,y_S$ which contains any
points added to $\Vis(R)$ between $\sigma+1$ and $\tau$.
\end{proof}

\noindent To finish the proof of lemma \ref{lem:cpl} just take an instance of
the coupling, and run it until $R$ makes $t$ moves. Since at each step
of the coupling $x_R = x_S$, $y_R=y_S$ and $z_R\leq z_S$ we get that
each time $S$ hit a specific floor point, $R$ hits it as well. The fact
that $z_R-\Wait_R = z_S - \Wait_S$ implies that in the remaining steps
$S$ has to complete $t$ moves, he can at most reach the same $z$
coordinate as $R$ (if he goes straight down), so he cannot bypass the
number of times $R$ hits $v$.
\end{proof}

We now use the coupling lemma for the following useful
corollary.

\begin{cor} \label{no_place}
For any floor point $v$ and any $t\in\N$,
$\E(V(t;v)\, |\, V(t;v)\neq 0) \leq \E V(t;\underline{0})$.
\end{cor}
\begin{proof}
Divide the probability space of all possible $t$-histories
according to the path the walk takes until reaching $v$ for the
first time (Or not reaching it at all). Examine one path which
does reach $v$ in $t$ steps, and let $T_v$ denote the first time
it hits $v$. Then, according to the coupling lemma, regardless of
the history until time $T_v$ ($\Vis(T_v)$), the expected number
of times the walk will hit $v$ in the \emph{next $t$ steps} is less
or equal to the expected number of times a walk starting at $v$
with no history (thus with no visited vertices) will visit $v$ in
$t$ steps. But clearly there is no difference between the (expected)
number of times a walk starting at $v$ will hit $v$ and the number
of times a walk starting at $\underline{0}$ hits $\underline{0}$,
so $\E(V(t;v)\, | \,T_v, \Vis(T_v)) \leq \E($number of returns
to $v$ in the next $t$ steps after $T_v \, | \, \Vis(T_v)) \leq
\E(V(t;\underline{0}))$. Since this holds for any value of
$T_v$ and $\Vis(T_v)$, the corollary follows.
\end{proof}

\section{Lower bound and proof of recurrence}

In the first few steps we will use the upper bound on the number of visits to a point
to get a lower bound on the number of new vertices the walk
visits, and consequently a lower bound on the number of times the
walk hits the floor. Here we still do not need the coupling argument.

\begin{lem}
Denote by $N(t)$ the number of different vertices the walk visits
till time $t$, then  $\exists c > 0$ such that $\E(N(t))
\geq ct/ \sqrt{\log t}$.
\end{lem}
\begin{proof}
Fix $t$ and then
\begin{align*} \E(N(t)) &= \E(\#v : V(v) \neq 0) = \sum_v\PP(V(v)\neq 0) =
\intertext{by the definition of conditional expectation,}
 &= \sum_v\frac{\E(V(v))}{\E(V(v) | V(v) \neq 0)} \geq
\intertext{by corollary \ref{cor:cond},}
&\geq c\sum_v \frac{\E(V(v))}{ \sqrt{\log t}} =
\frac{c}{ \sqrt{\log t}}\sum_v\E(V(v)) = \\
&=\frac{c}{\sqrt{\log t}}\E\Big(\sum_v V(v)\Big) =
 \frac{ct}{\sqrt{\log t}}.\qedhere
\end{align*}
\end{proof}

We denote by $DF(t)$ the number of \emph{different floor} points
the walk visits till time $t$. The next lemma bounds $\E(DF(t))$.
\begin{lem}\label{lem:DFn}
$\E(DF(t)) < Ct/\log t$
\end{lem}
\begin{proof}
 Since any two distinct floor points have different $(x,y)$
 coordinates, the number of different floor points visited by the
 ERW (till time $t$), is bounded from above by the number of different points
 its $(x,y)$ projection visits. But the projection of the ERW
 on the $(x,y)$ plane is a simple random walk of length $\leq
 t$. Therefore  $\E(DF(t))$ is bounded above by the expected number of
 different vertices a SRW visits in $t$ steps, which by
Dvoretzky-Erd\"os \cite{DE51} is $< Ct/\log t$.
\end{proof}

\begin{cor}
There exists a 
$c>0$ such that $\E(|\Vis(t)|) \geq ct/\sqrt{\log t}$.
\end{cor}
\begin{proof}
Because $|\Vis(t)|=N(t)-DF(t)$.
\end{proof}

Now we are ready to bound from below the expected number of times
the ERW hits the floor.
\begin{lem}\label{lem:Fn}
Denote by $F(t)$ the number of times the ERW hits the floor
($z=0$) until time $t$. Then there is a positive constant $c>0$,
independent of $t$ such that $\E(F(t)) \geq ct/\sqrt{\log
t}$.
\end{lem}
\begin{proof}
Look at the expected change of the $z$ coordinate when the walk
makes a single step. The walk has one of three behaviors according
to the type of vertex it is currently in. If the walk is in a
visited vertex --- it acts as a SRW, so
the expected change to the $z$ coordinate is $0$.
If the walk is in a new vertex --- it goes down a step ---
so the expected change to the $z$ coordinate is $-1$, and finally,
if the walk is on the floor then the expected change in the $z$
coordinate is $\frac{1}{5}$. We get, by linearity of expectation,
\[
\E(z(t)) = -1 \cdot \E(|\Vis(t)|) + \frac{1}{5}\cdot\E(F(t)).
\]
Since the walk always stays on the upper half space $z\geq 0$, we have
$\E(z(t)) \geq 0$ so $\E(F(t)) \geq 5\E(|\Vis(t)|) \geq 5ct/\sqrt{\log
t}$ by the previous corollary.
\end{proof}

With the above estimates, all we need to do is combine our bounds on $\E(F(t))$ and
$\E(DF(t))$ with corollary \ref{no_place} to
get:

\begin{thm}\label{thm:lower}
There exists a constant $c>0$ such that  $\E(V(t;\underline{0}))
\geq c\sqrt{\log t}$.
\end{thm}
\begin{proof}
From linearity of expectation,
\[
\E(F(t)) = \sum_{v\in \text{floor}}\E(V(v)) = \sum_{v\in
\text{floor}}(\E(V(v) \, | \, V(v)\neq 0) \cdot \PP(V(v) \neq 0).
\]
By corollary \ref{no_place}, $\E(V(v) \,|\, V(v)\neq 0) \leq \E
V(\underline{0})$ for any floor point $v$, so
\begin{align*}
\E(F(t)) &\leq \sum_{v\in \text{floor}}\E(V(\underline{0})) \cdot
\PP(V(v) > 0) = \E(V(\underline{0})) \cdot \sum_{v\in
\text{floor}}\PP(V(v)>0) = \\
 &= \E(V(\underline{0})) \cdot \E(DF(t))
\end{align*}
so by lemmas \ref{lem:DFn} and \ref{lem:Fn},
\[
\E(V(\underline{0})) \geq \frac{\E(F(t))}{\E(DF(t))} \geq
\frac{c t/\sqrt{\log t}}{C t/\log t} = c \sqrt{\log t}
\]
For some  absolute constant $c>0$.
\end{proof}

\begin{thm}
ERW is recurrent.
\end{thm}
\begin{proof}
Assume to the contrary that there is a positive probability that
an ERW $R$ visits the origin exactly $k$ times, $k$ finite. Let $\tau$
be the hitting time when $R$ reaches the origin for the $k$-th time.
Examining $\Vis_R(\tau)$ we see that there is a finite set $X$ and a
positive probability $p_0>0$ such that if the walk $R$ is
currently at $\underline{0}$ and $\Vis_R=X$ then $\PP(R
\text{ will not return to }\underline{0}) \geq p_0$.
By the coupling lemma (lemma \ref{lem:cpl}), if the walk $R$ is at
$\underline{0}$, and $X\subset \Vis(R)$, then
$\PP(R\text{ will not return again to }\underline{0}) \geq p_0$.
Since the set of visited vertices of a walk only increases, we
conclude that once $X\subset\Vis_R$, the number of times
the walk returns to $\underline{0}$ is dominated by a geometric
random variable with parameter $1-p_0$.

Next, notice that for any finite set of vertices $Y$, there is a
positive probability $p_1(Y)$ such that a walk currently at
$\underline{0}$ will visit all the vertices in $Y$ before
returning to $\underline{0}$ with probability $\geq p_1$,
regardless of the walk's history. (One such possible path is simply
reaching each vertex $v=(x,y,z)\in Y$ by walking on the floor
[avoiding $(0,0)$] till you reach $(x,y,0)$, then climbing slowly
up until reaching the desired vertex, and finish the path by
dropping on the $(0,0,*)$ column from high enough.)
Thus the number of returns of ERW to $\underline{0}$ before all vertices
of $X$ are visited is dominated by a geometric random variable
with parameter $1-p_1(X)$.

Combining the above we get that the total number of returns of ERW
to $\underline{0}$ is dominated by the sum of two geometric variables
(with parameters $1-p_0,1-p_1)$, and thus has finite expectation. But this
contradicts theorem \ref{thm:lower}, (The expected number of returns till time
$t$ behaves like $\sqrt{\log t}$), so ERW is recurrent.
\end{proof}

\end{document}